\documentstyle[11pt,twoside]{article}
\include{graphicx}
\title{A note on an integral associated with the Kelvin ship-wave pattern}
\author{\sc R. B.\ Paris \\
{\em Division of Computing and Mathematics}, \\
{\em University of Abertay Dundee, Dundee DD1 1HG, UK}
}
\begin{document}
\def\f#1#2{\mbox{${\textstyle \frac{#1}{#2}}$}}
\def\dfrac#1#2{\displaystyle{\frac{#1}{#2}}}
\def\boldal{\mbox{\boldmath $\alpha$}}
{\newcommand{\Sgoth}{S\;\!\!\!\!\!/}
\newcommand{\bee}{\begin{equation}}
\newcommand{\ee}{\end{equation}}
\newcommand{\lam}{\lambda}
\newcommand{\ka}{\kappa}
\newcommand{\al}{\alpha}
\newcommand{\fr}{\frac{1}{2}}
\newcommand{\fs}{\f{1}{2}}
\newcommand{\g}{\Gamma}
\newcommand{\br}{\biggr}
\newcommand{\bl}{\biggl}
\newcommand{\ra}{\rightarrow}
\newcommand{\mbint}{\frac{1}{2\pi i}\int_{c-\infty i}^{c+\infty i}}
\newcommand{\mbcint}{\frac{1}{2\pi i}\int_C}
\newcommand{\mboint}{\frac{1}{2\pi i}\int_{-\infty i}^{\infty i}}
\newcommand{\gtwid}{\raisebox{-.8ex}{\mbox{$\stackrel{\textstyle >}{\sim}$}}}
\newcommand{\ltwid}{\raisebox{-.8ex}{\mbox{$\stackrel{\textstyle <}{\sim}$}}}
\renewcommand{\topfraction}{0.9}
\renewcommand{\bottomfraction}{0.9}
\renewcommand{\textfraction}{0.05}
\newcommand{\mcol}{\multicolumn}
\date{}
\maketitle
\pagestyle{myheadings}
\markboth{\hfill \sc R. B.\ Paris  \hfill}
{\hfill \sc  Kelvin ship-wave pattern\hfill}
\begin{abstract}
The velocity potential in the Kelvin ship-wave source can be partly expressed in terms of space derivatives of
the single integral
\[F(x,\rho,\alpha)=\int_{-\infty}^\infty \exp\,[-\fs\rho \cosh (2u-i\alpha)] \cos (x\cosh u)\,du,\]
where $(x, \rho, \alpha)$ are cylindrical polar coordinates with origin based at the source and $-\fs\pi\leq\alpha\leq\fs\pi$. An asymptotic expansion of $F(x,\rho,\alpha)$ when $x$ and $\rho$ are small, but such that $M\equiv x^2/(4\rho)$ is large, was given using a non-rigorous approach by Bessho in 1964 as a sum involving products of Bessel functions. This expansion, together with an additional integral term, was subsequently proved by Ursell in 1988. 

Our aim here is to present an alternative asymptotic procedure for the case of large $M$. The resulting expansion
consists of three distinct parts: a convergent sum involving the Struve functions, an asymptotic series and an exponentially small saddle-point contribution. Numerical computations are carried out to verify the accuracy of our expansion.
\vspace{0.4cm}

\noindent {\bf Mathematics Subject Classification:} 30E15, 34E05, 41A30, 41A60, 76B20 
\vspace{0.3cm}

\noindent {\bf Keywords:} Kelvin ship-wave pattern, asymptotic expansion,  Struve functions
\end{abstract}

\vspace{0.3cm}

\noindent $\,$\hrulefill $\,$

\vspace{0.2cm}

\begin{center}
{\bf 1. \  Introduction}
\end{center}
\setcounter{section}{1}
\setcounter{equation}{0}
\renewcommand{\theequation}{\arabic{section}.\arabic{equation}}
The velocity potential in the Kelvin ship-wave source can be partly expressed in terms of space derivatives of
the single integral
\[\int_{-\infty}^\infty \exp\,[-\fs z\cosh 2u+\fs iy \sinh 2u] \cos (x\cosh u)\,du\hspace{4cm}\]
\bee\label{e11}
\hspace{4cm}=\int_{-\infty}^\infty \exp\,[-\fs\rho \cosh (2u-i\alpha)] \cos (x\cosh u)\,du\equiv F(x, \rho, \alpha),
\ee
where the origin of cylindrical polar $(\rho,\alpha,x)$ coordinates is at the source, $x>0$ is the non-dimensional horizontal coordinate along the track, $y=\rho \sin \alpha$ is the horizontal coordinate perpendicular to the track and $z=\rho \cos \alpha$ is the vertical coordinate (increasing with depth). The polar angle $\alpha$ is zero
directly below the origin and satisfies $-\fs\pi\leq\alpha\leq\fs\pi$, although it is sufficient to restrict attention to $0\leq\alpha\leq\fs\pi$ since $F(x, \rho, \alpha)$ is an even function of $\alpha$.

The integral $F(x, \rho, \alpha)$ is difficult to evaluate when $x$ and $\rho$ are small. In 1964, Bessho \cite{B} gave the expansion
\[F(x,\rho,\alpha)=K_0(\fs\rho) J_0(x)+2\sum_{m=1}^\infty (-)^m \cos m\alpha\,K_m(\fs\rho) J_{2m}(x),\]
which he obtained by expanding the cosine factor as a series of Bessel functions. Here $J$ and $K$ denote the usual Bessel functions, and from the well-known large-order behaviour
\[J_{2m}(x)\sim \frac{1}{2\sqrt{\pi m}} \bl(\frac{ex}{4m}\br)^m,\qquad K_m(\fs\rho)\sim \sqrt{\frac{\pi}{2m}} \bl(\frac{e\rho}{4m}\br)^{-m}\quad(m\ra +\infty)\]
it follows that the late terms in the above sum are controlled by $(ex^2/4\rho)^m/m^{m+1}$. The expansion is therefore absolutely convergent, but in the case where $x$ and $\rho$ are small with $x^2/(4\rho)$ assumed to be large, the convergence is slow and the accuracy is reduced by the accumulation of cancellation errors.

An expansion suitable for $M\equiv x^2/(4\rho)\gg 1$ was established rigorously by Ursell \cite{U1} in the form
\[F(x, \rho, \alpha)=-\pi \bl(I_0(\fs\rho) Y_0(x)+2\!\sum_{m\leq M} \cos m\alpha\,I_m(\fs\rho) Y_{2m}(x)\br) \]
\bee\label{e12}
+\Re\,\int_{-\infty+\frac{1}{2}i|\alpha|}^{-\infty+\frac{1}{2}i|\alpha|-\pi i} \exp\,[-\fs\rho \cosh (2u-i|\alpha|)+ix \cosh u]\,du +O(e^{-M}),
\ee
where $I$ and $Y$ are the usual Bessel functions. The series extended over infinite values of $m$ (without the integral and order term) had been obtained in a non-rigorous fashion in \cite{B}. From the large-$m$ behaviour
\[I_m(\fs\rho)\sim \frac{1}{\sqrt{2\pi m}}\bl(\frac{e\rho}{4m}\br)^m   ,\qquad Y_{2m}(x)\sim -\frac{1}{\sqrt{\pi m}}\bl(\frac{ex}{4m}\br)^{-2m}     \quad (m\ra+\infty),\]
it is easily seen that the terms in the series in (\ref{e12}) are comparable to an expansion in inverse powers of $M$. The additional term represented by the integral in (\ref{e12}) can be estimated by the method of steepest descents and is given by \cite[Eq.~(2.26)]{U1}
\bee\label{e13}
(\pi/M)^\fr \,e^{-(M-\fr\rho) \cos \alpha}\,\sin\,[(M+\fs\rho) \sin \alpha+\fs\alpha],
\ee
which is exponentially small except when $\alpha$ is near $\fs\pi$. An algorithm for a more precise estimation of this term is given in \cite{N}.

Our aim in this note is to present a different asymptotic procedure for $F(x. \rho, \alpha)$ when $x$ and $\rho$ small, but with
the parameter $M\gg 1$. We deal first with the expansion in the midplane (that is when $\alpha=0$) and then extend the procedure for $0<\alpha\leq\fs\pi$. The expansion we obtain consists of three parts: the first is a convergent series involving  Struve functions that is rapidly convergent for small values of $x$ and $\rho$; the second is an asymptotic series and the third is an exponentially small saddle-point contribution. This latter contribution agrees with (\ref{e13}) when $\alpha$ is bounded away from zero. Numerical calculations reveal that our expansion and that in (\ref{e12}) produce comparable accuracy.

\vspace{0.6cm}

\begin{center}
{\bf 2. \ Asymptotic evaluation when $\alpha=0$}
\end{center}
\setcounter{section}{2}
\setcounter{equation}{0}
\renewcommand{\theequation}{\arabic{section}.\arabic{equation}}
We consider the case where $x$ and $\rho$ are small but $x^2/(4\rho)$ is assumed to be large.
When $\alpha=0$, we have from (\ref{e11})
\bee\label{e21}
F(x, \rho, 0)=\int_{-\infty}^\infty e^{-\fr\rho \cosh 2u}\cos (x\cosh u)\,du
=2e^{\fr\rho} \int_p^\infty \frac{e^{-M\tau^2}\cos\,(2M\tau)}{\sqrt{\tau^2-p^2}}\,d\tau,
\ee
where we have set $\cosh u=x\tau/(2\rho)$ and defined
\[M=\frac{x^2}{4\rho}~,\qquad p=\frac{2\rho}{x}=\frac{x}{2M}.\]

The integrals
\bee\label{e22}
\int_p^\infty \frac{e^{-M(\tau^2\mp 2i\tau)}}{\sqrt{\tau^2-p^2}}\,d\tau
\ee
have saddle points at $\tau=\pm i$, with the steepest descent paths through these points being the doubly infinite lines parallel to the real $\tau$-axis. The parts of these paths situated in $\Re (\tau)\geq 0$ are labelled $L_1$ and $L_2$ in Fig.~1. The square root in the denominator makes it necessary to introduce a branch cut the $\tau$-plane along $[-p, p\,]$. 

Considering the integral with the upper sign in (\ref{e22}), we deform the integration path $[p, \infty)$ into a semi-circular path round the branch point $\tau=p$, the upper side of the branch cut to the origin followed by the positive imaginary axis $0\leq \Im\,(\tau)\leq 1$ and then out to infinity along the steepest descent path (which we label $L_1$) through $\tau=i$. This corresponds to the paths shown in the upper half of Fig.\,1. 
If the phase of $(\tau^2-p^2)^{1/2}$ is zero on the positive axis on the right of $\tau=p$, then the phase is $\pm\fs\pi$ on the upper and lower sides of the branch cut, and on the positive and negative halves of the imaginary axis, respectively. Then we find, on letting the radius of the indentation around $\tau=p$ shrink to zero,
\[\int_p^\infty \frac{e^{-M(\tau^2-2i\tau)}}{\sqrt{\tau^2-p^2}}\,d\tau
=i\int_0^p\frac{e^{-M(\tau^2-2i\tau)}}{\sqrt{p^2-\tau^2}}\,d\tau+
\int_0^1 \frac{e^{M(u^2-2u)}}{\sqrt{u^2+p^2}}\,du
-i\int_{L_1}\frac{e^{-M(\tau^2-2i\tau)}}{\sqrt{p^2-\tau^2}}\,d\tau
\]
Similar considerations applied to the integral with the lower sign leads to the above relation with $i$ replaced by $-i$ and the path $L_1$ by $L_2$. 
Then we find
\bee\label{e23}
e^{-\fr\rho} F(x, \rho, 0)=-2\int_0^p \frac{e^{-M\tau^2} \sin\,(2M\tau)}{\sqrt{p^2-\tau^2}}\,d\tau+2\int_0^1
\frac{e^{M(u^2-2u)}}{\sqrt{u^2+p^2}}\,du+I_s,
\ee
where $I_s$ denotes the saddle-point contribution
\bee\label{e23a}
I_s=i\int_{L_2} \frac{e^{-M(\tau^2+2i\tau)}}{\sqrt{p^2-\tau^2}}\,d\tau-i\int_{L_1} \frac{e^{-M(\tau^2-2i\tau)}}{\sqrt{p^2-\tau^2}}\,d\tau.
\ee
\begin{figure}[t]
	\begin{center}\includegraphics[width=0.4\textwidth]{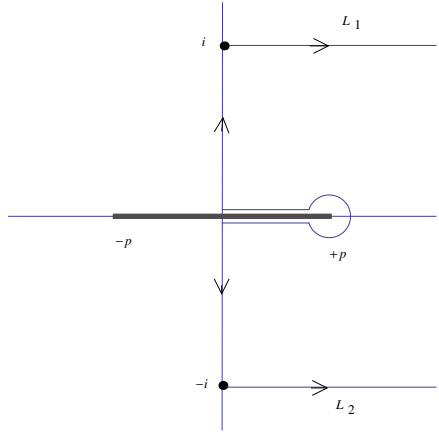}\\
\caption{\small{The deformation of the integration path $[p, \infty)$ in the $\tau$-plane through the saddles at $\pm i$. There is a branch cut along $[-p,p\,]$.}}
	\end{center}
\end{figure}

We now consider in turn each of the integral contributions in (\ref{e23}). For this we shall require the evaluation\footnote{This result can be obtained by using the series expansion of the exponential and evaluating the resulting integrals in terms of the Beta function. The result given in \cite[3.478(3)]{GR} incorrectly excludes this case.}
\bee\label{e24}
\int_0^p e^{-M\tau^2} \tau^{2k+1} (p^2-\tau^2)^{\mu-1}d\tau=\frac{p^{2k+2\mu}k! \g(\mu)}{2\g(k+\mu+1)}{}_1F_1(k+1;k+\mu+1;-Mp^2)
\ee
for $\mu>0$ and $k=0, 1, 2, \ldots $, where ${}_1F_1$ denotes the confluent hypergeometric function.
\vspace{0.4cm}

\noindent {2.1\it \ \ The branch-cut contribution}
\vspace{0.2cm}

\noindent We first consider the branch-cut contribution
\[I_1=\int_0^p \frac{e^{-M\tau^2} \sin\,(2M\tau)}{\sqrt{p^2-\tau^2}}\,d\tau.\]
Expanding the sine function followed by reversal of the order of summation and integration, we find
\begin{eqnarray*}
I_1&=&\sum_{k=0}^\infty\frac{(-)^k (2M)^{2k+1}}{(2k+1)!}\int_0^p\frac{e^{-M\tau^2}\tau^{2k+1}}{\sqrt{p^2-\tau^2}}\,d\tau\\
&=&\frac{\sqrt{\pi}}{2}\sum_{k=0}^\infty \frac{(-)^k x^{2k+1}}{(2k+1)!}
\,\frac{k!}{\g(k+\f{3}{2})}\,{}_1F_1(k+1;k+\f{3}{2};-\rho)\\
&=&\frac{\pi e^{-\rho}}{2}\sum_{k=0}^\infty \frac{(-)^k (x/2)^{2k+1}}{\g^2(k+\f{3}{2})}\,{}_1F_1(\fs; k+\f{3}{2};\rho)
\end{eqnarray*}
by (\ref{e24}) with $\mu=\fs$, $Mp^2=\rho$ and $2Mp=x$. Here we have employed Kummer's transformation for the confluent hypergeometric function and also the duplication formula for the gamma function.

If the ${}_1F_1$ function is expressed in its series form, we then obtain
\begin{eqnarray}
I_1&=&\frac{\pi e^{-\rho}}{2}\sum_{r=0}^\infty\frac{(\fs)_r}{r!}\,\rho^r \sum_{k=0}^\infty \frac{(-)^k (\fs x)^{2k+1}}{\g(k+\f{3}{2}) \g(k+r+\f{3}{2})}\nonumber\\
&=&\frac{\pi e^{-\rho}}{2}\sum_{r=0}^\infty\frac{(\fs)_r}{r!}\rho^r {\cal H}_r(x),\label{e25}
\end{eqnarray}
where $(a)_r=\g(a+r)/\g(a)$ is Pochhammer's symbol and ${\cal H}_\nu(x)$ is the `scaled' Struve function defined 
by \cite[p.~288]{DLMF}
\bee\label{e2s}
{\cal H}_\nu(x):=(\fs x)^{-\nu} {\bf H}_\nu(x)=\sum_{k=0}^\infty \frac{(-)^k (\fs x)^{2k+1}}{\g(k+\f{3}{2})\g(k+\nu+\f{3}{2})}.
\ee
The series on the right-hand side of (\ref{e25}) is absolutely convergent, since 
the behaviour of the scaled Struve function as $r\ra \infty$ is given by  ${\cal H}_r(x)\sim xe^r/(\pi \surd 2\,r^{r+1})$; see \cite[Eq.~(11.6.5)]{DLMF}.
\vspace{0.4cm}

\noindent {2.2\it\ \ The imaginary axis contribution}
\vspace{0.2cm}

\noindent We now consider the second integral in (\ref{e23}), namely
\[I_2=\int_0^1\frac{e^{M(u^2-2u)}}{\sqrt{u^2+p^2}}\,du.\]
This represents the contribution from the origin down the steepest descent path (relative to the origin) to the saddles at $\tau=\pm i$ and is a Laplace integral with a linear endpoint at $u=0$. However, it is not possible to employ the standard asymptotic estimate for such an integral as described in, for example, \cite[2.4(iii)]{DLMF}, \cite[p.~14]{P}, since we have $M\ra\infty$ and $p=x/(2M)\ra 0$.
With the change of variable $\xi=u/p$, this becomes
\[I_2=\int_0^{\xi_0} \frac{e^{-x\xi+\rho\xi^2}}{\sqrt{1+\xi^2}}\,d\xi, \qquad \xi_0=\frac{1}{p}=\frac{2M}{x}.\]

We now write $e^{\rho\xi^2}$ as a finite series of $n$ ($>1$) terms together with a remainder expressed in Lagrange's form \cite[p.~328]{C}, \cite[p.~96]{WW}:
\bee\label{e2L}
e^{\rho\xi^2}=\sum_{k=0}^{n-1}\frac{(\rho\xi^2)^k}{k!}+\frac{(\rho\xi^2)^n}{n!}\,e^{\theta\rho\xi^2}\qquad (0<\theta<1).
\ee
This then enables us to write
\bee\label{e26}
I_2=\sum_{k=0}^{n-1}\frac{\rho^k}{k!} \bl(\int_0^\infty-\int_{\xi_0}^\infty\br) \frac{\xi^{2k} e^{-x\xi}}{\sqrt{1+\xi^2}}\,d\xi+R_n,
\ee
where
\bee\label{e26a}
R_n=\frac{\rho^n}{n!}\int_0^{\xi_0} \frac{\xi^{2n} e^{-x\xi+\theta\rho\xi^2}}{\sqrt{1+\xi^2}}\,d\xi.
\ee

The main contribution to the finite sum in (\ref{e26}) may be written in the form
\bee\label{e29}
S=\frac{\pi}{2}\sum_{k=0}^{n-1}\frac{M^{-k}}{2^{2k} k!}\,C_k(x),\qquad C_k(x):=\frac{2}{\pi}\,x^{2k} \!\!\int_0^\infty\frac{\xi^{2k}e^{-x\xi}}{\sqrt{1+\xi^2}}\,d\xi.
\ee
To evaluate the integrals $C_k(x)$, we observe that for integer $m\geq 0$ \cite[Eq.~(11.5.2)]{DLMF}
\[\sum_{r=0}^m \bl(\!\!\begin{array}{c}m\\r\end{array}\!\!\br) \int_0^\infty \frac{\xi^{2r} e^{-x\xi}}{\sqrt{1+\xi^2}}\,d\xi=\int_0^\infty (1+\xi^2)^{m-\fr} e^{-x\xi}\,d\xi=\frac{\surd\pi\g(m+\fs)}{2^{1-m}}\, x^{-m}\,{\bf K}_m(x),\] 
\[{\bf K}_\nu(z)={\bf H}_\nu(z)-Y_\nu(z),\] 
where ${\bf K}_\nu(x)$ and ${\bf H}_\nu(x)$ are Struve functions and $Y_\nu(z)$ the usual Bessel function; see \cite[Eq.~(11.2.5)]{DLMF}. Then, with ${\cal K}_m(x)\equiv x^m {\bf K}_m(x)$, we obtain the recurrence relation
\bee\label{e210}
C_m(x)=2^m(\fs)_m\,{\cal K}_m(x)-\sum_{r=0}^{m-1} \bl(\!\!\begin{array}{c}m\\r\end{array}\!\!\br)\,x^{2(m-r)}C_r(x)\qquad (m\geq 1),
\ee
from which we find the first few $C_k(x)$ given by
\[C_0(x)={\cal K}_0(x),\quad C_1(x)={\cal K}_1(x)-x^2{\cal K}_0(x),\]
\[C_2(x)=3{\cal K}_2(x)-2x^2{\cal K}_1(x)+x^4{\cal K}_0(x),\]
\[C_3(x)=15{\cal K}_3(x)-9x^2{\cal K}_2(x)+3x^4{\cal K}_1(x)-x^6{\cal K}_0(x),\]
\[C_4(x)=105{\cal K}_4(x)-60x^2{\cal K}_3(x)+18x^4{\cal K}_2(x)-4x^6{\cal K}_1(x)+x^8{\cal K}_0(x)~.\]

In Appendix A it is established that the contribution resulting from the integral over $[\xi_0, \infty)$ in (\ref{e26}) is O($e^{-M})$, when\footnote{If $n$ is finite, this estimate can be replaced by $O(M^n e^{-2M})$; see (\ref{a4}).} $n<M$, and that $R_n=O(M^{-n})$. Then, from (\ref{e29}) we have the result
\bee\label{e211}
I_2=\frac{\pi}{2}\bl\{\sum_{k=0}^{n-1}\frac{M^{-k}}{2^{2k} k!}\,C_k(x)+O(M^{-n})\br\}+O(e^{-M})
\qquad (n<M)
\ee
as $M\ra\infty$.
\vspace{0.4cm}

\noindent 2.3 {\it\ \  The saddle-point contribution}
\vspace{0.2cm}

\noindent The final term to consider in (\ref{e23}) is the saddle-point contribution $I_s$ defined in (\ref{e23a}). Recalling that the integration paths $L_1$ and $L_2$ are $\tau=u\pm i$, $0\leq u\leq\infty$, we have
\[I_s=ie^{-M} \int_0^\infty e^{-Mu^2} f(u)\,du,\qquad f(u)=\frac{1}{\sqrt{p^2-(u-i)^2}}-\frac{1}{\sqrt{p^2-(u+i)^2}}.\]
%where
%\begin{eqnarray*}
%f(u)&=&\frac{1}{\sqrt{p^2-(u-i)^2}}-\frac{1}{\sqrt{p^2-(u+i)^2}}\\
%&=&-\frac{2iu}{(1+p^2)^{3/2}}+O(u^3).
%\end{eqnarray*}
Then, since $f(u)=2iu/(1+p^2)^{3/2}+O(u^3)$ as $u\ra 0$, we obtain
\bee\label{e212}
I_s=\frac{2e^{-M}}{(1+p^2)^{3/2}} \int_0^\infty e^{-Mu^2}\{u+O(u^3)\}\,du=\frac{e^{-M}}{M(1+p^2)^{3/2}}\{1+O(M^{-1})\}
\ee
as $M\ra\infty$.
\bigskip

From (\ref{e23}), (\ref{e25}), (\ref{e211}) and (\ref{e212}), we obtain our desired result when $\alpha=0$ given
in the following theorem.
\newtheorem{theorem}{Theorem}
\begin{theorem}$\!\!\!.$\ \ Let $x$ and $\rho$ be positive and bounded (to fix ideas suppose $0<x<1$ and $0<\rho<1$) and let $\alpha=0$. Then, when $M\equiv x^2/(4\rho)$ is large we have
\[F(x,\rho,0)=-\pi e^{-\fr\rho}\sum_{r=0}^\infty \frac{(\fs)_r}{r!}\,\rho^r {\cal H}_r(x)+
\pi e^{\fr\rho}\bl\{\sum_{k=0}^{n-1} \frac{M^{-k}}{2^{2k} k!}\,C_k(x)+O(M^{-n})\br\}\]
\bee\label{e213}
\hspace{4cm}+O(e^{-M})+O(M^{-1}e^{-M})
\ee 
where $n<M$ is a fixed positive integer, ${\cal H}_r(x)$ is a scaled Struve function and the coefficients $C_k(x)$ are defined by the recurrence in (\ref{e29}).
\end{theorem}
\vspace{0.6cm}

\begin{center}
{\bf 3. \ Asymptotic evaluation when $0<\alpha\leq\fs\pi$}
\end{center}
\setcounter{section}{3}
\setcounter{equation}{0}
\renewcommand{\theequation}{\arabic{section}.\arabic{equation}}
When $0<\alpha\leq\fs\pi$, we can make the substitution $u \ra u+\fs i\alpha$ in (\ref{e11}), appealing to Cauchy's theorem to place the integration path back onto $(-\infty, \infty)$, to find
\[F(x,\rho,\alpha)=\int_{-\infty}^\infty e^{-\fr\rho \cosh 2u} \cos (xc \cosh u+isx \sinh u)\,du,\]
where, for brevity in what follows, we have set
\[c\equiv \cos \fs\alpha,\qquad s\equiv \sin \fs\alpha.\]
With $\cosh u=x\tau/(2\rho)$, some routine algebra then shows that
\[e^{-\fr\rho} F(x,\rho,\alpha)=J_+(M, \alpha)+J_-(M,\alpha),\]
where
\bee\label{e31}
J_\pm(M,\alpha):=\int_p^\infty e^{M(-\tau^2\pm 2ic\tau)}\ \frac{\cosh(2Ms\sqrt{\tau^2-p^2})}{\sqrt{\tau^2-p^2}}\,d\tau.
\ee

The integral $J_+(M,\alpha)$ contains two components associated with the exponential factors $\exp [M\psi(\tau)]$, where
\[\psi(\tau)=-\tau^2+2ic\tau\pm 2s\sqrt{\tau^2-p^2};\]
we shall refer to these as the `positive' and `negative' components of $J_+(M,\alpha)$, respectively. The components are associated with saddle points given by
\bee\label{e35a}
\tau_s=ic\pm s \sqrt{\frac{\tau_s^2}{\tau_s^2-p^2}}.
\ee
Remembering that $p=x/(2M)$, we see that for large $M$ each component of $J_+(M,\alpha)$ has a saddle situated at
\bee\label{e35b}
\tau_s^\pm\simeq ie^{\mp i\alpha/2}-\fs p^2s e^{\pm i\alpha} \qquad (p\ra 0).
\ee
Some routine algebra, using (\ref{e35a}) and (\ref{e35b}), shows that
\begin{eqnarray}
\psi(\tau_s^\pm)&=&-\tau_s^\pm{}^2+2ic\tau_s^\pm +\frac{2s^2\tau_s^\pm}{\tau_s^\pm-ic}\simeq -\tau_s^\pm{}^2+2ie^{\mp i\alpha/2}\pm ip^2s e^{\pm i\alpha/2}\nonumber\\
&\simeq& -e^{\mp i\alpha}\pm ip^2 se^{\pm i\alpha/2} \qquad (p\ra 0).\label{e35e}
\end{eqnarray}

\begin{figure}[t]
	\begin{center}
	{\tiny($a$)}\ \includegraphics[width=0.25\textwidth]{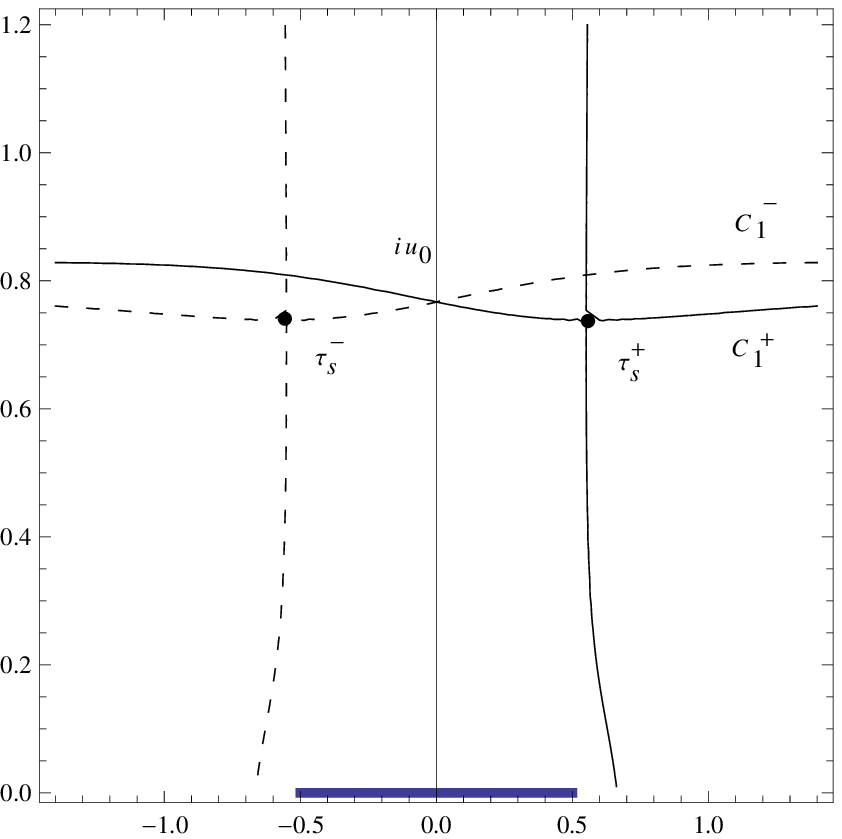}\qquad {\tiny($b$)} \includegraphics[width=0.25\textwidth]{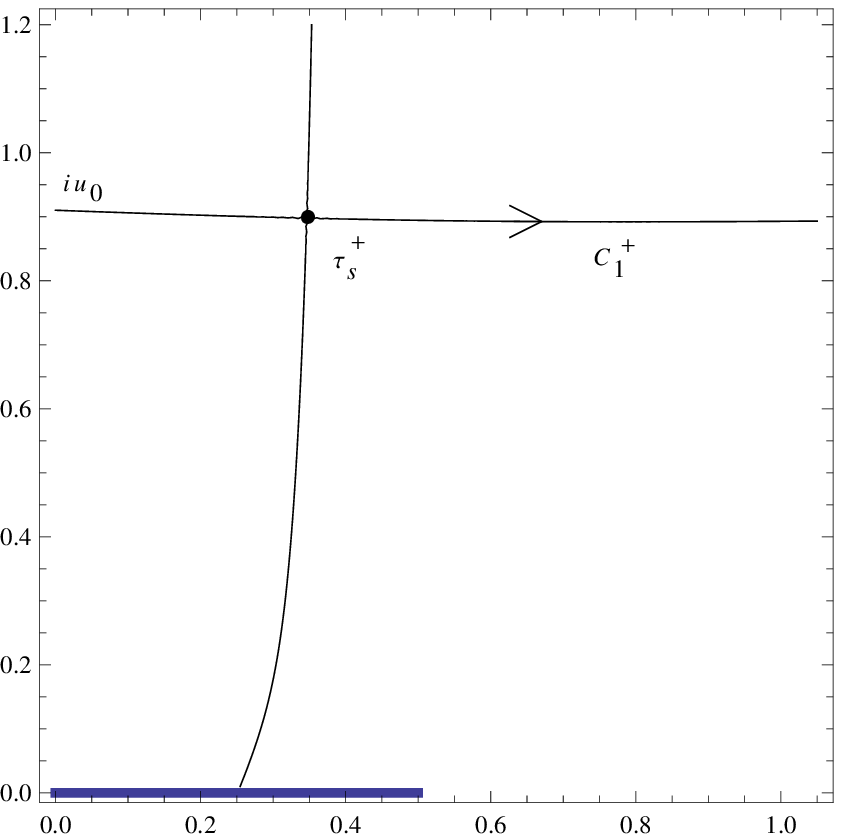} \qquad 
{\tiny($c$)}\includegraphics[width=0.25\textwidth]{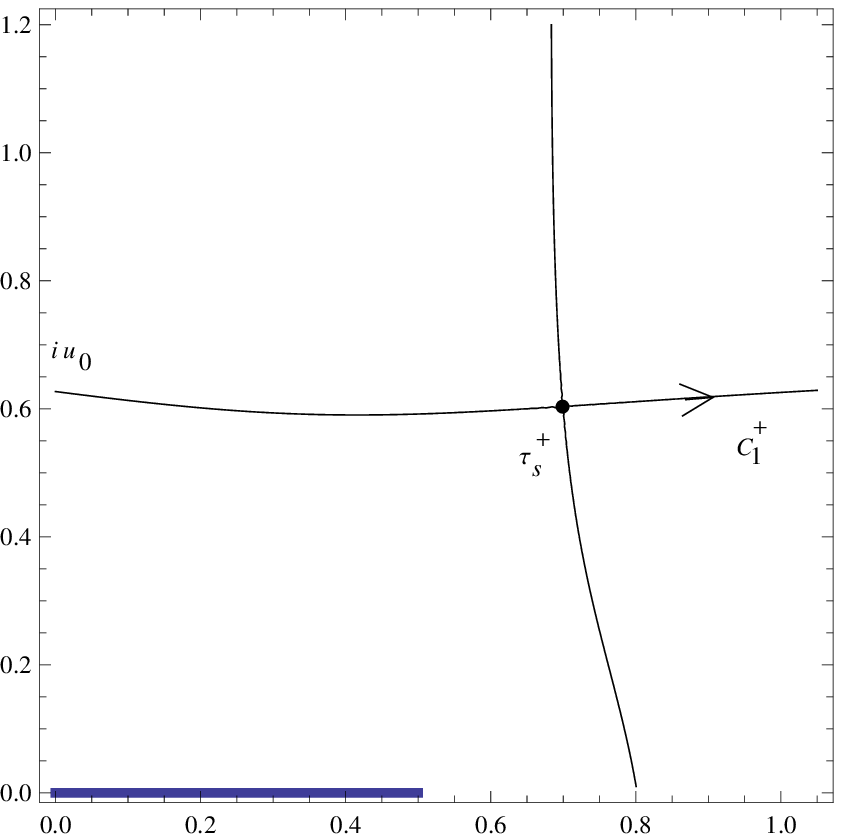}\\
\caption{\small{(a) The steepest descent (${\cal C}_1^\pm$) and ascent paths through the saddles $\tau_s^\pm$ for the two components of $J_+(M,\alpha)$ and the point $iu_0$ when $p=\fs$ and $\alpha=2\pi/5$. The steepest paths through $\tau_s^+$ in $\Re (\tau)\geq 0$ when (b) $\alpha=\pi/4$ and (c) $\alpha=\pi/2$. There is a branch cut along $[-p,p\,]$.}}
	\end{center}
\end{figure}

Typical paths of steepest descent through these two saddles are shown in Fig.~2(a), where that through the saddle  $\tau_s^-$ is shown dashed; only the saddle $\tau_s^+$ will be of relevance in our considerations. Both these paths, which we label ${\cal C}_1^+$ and ${\cal C}_1^-$ respectively, intersect the imaginary $\tau$-axis at the common point $\tau=iu_0$, where 
\[\Im (\psi(\tau_s^\pm))=\pm 2s \sqrt{u_0^2+p^2}\simeq\pm (\sin \alpha+\fs p^2 sc);\]
that is
\bee\label{e35}
\hspace{2cm} u_0\simeq c(1-\fs p^2 \tan^2 \fs\alpha) \qquad (0\leq\alpha\leq\fs\pi,\ \ p\ra 0).
\ee
The integral $J_-(M,\alpha)$ is associated with a conjugate set of saddles in the lower half-plane and conjugate steepest descent paths ${\cal C}_2^+$ and ${\cal C}_2^-$. When $\alpha=0$ the steepest descent paths become the horizontal lines $L_1$ and $L_2$ emanating from the points $\tau=\pm i$ shown in Fig.~1.

We now deform the integration path as described in Section 2, recalling that the phase of $(\tau^2-p^2)^{1/2}$ is $\pm\fs\pi$ on the upper and lower sides of the branch cut $[-p,p]$. The integral $J_+(M,\alpha)$ is taken along the upper side of the cut to the origin, followed by part of the positive imaginary axis $[0, iu_0]$ and then each component is taken out to infinity along its associated steepest descent path ${\cal C}_1^+$ and ${\cal C}_1^-$.
The positive component of $J_+(M,\alpha)$ passes over the saddle at $\tau_s^+$, whereas the negative component passes over part of the steepest descent path situated in $\Re (\tau)\geq 0$ that emanates from the saddle $\tau_s^-$ in $\Re (\tau)<0$.
The integral $J_-(M,\alpha)$ is taken along a similar path in the lower half-plane passing out to infinity along the paths ${\cal C}_2^+$ and ${\cal C}_2^-$. 
Then we find
\[e^{-\fr\rho} F(x,\rho,\alpha)=-2\int_0^p e^{-M\tau^2} \sin (2Mc\tau)\,\frac{\cos (2Ms\sqrt{p^2-\tau^2})}{\sqrt{p^2-\tau^2}}\,d\tau
\]
\bee\label{e32}
+2\int_0^{u_0} e^{M(u^2-2cu)}\,\frac{\cos (2Ms\sqrt{u^2+p^2})}{\sqrt{u^2+p^2}}\,du+I_s,
\ee
where
\[I_s=2\Im \!\int_{{\cal L}}e^{-M(\tau^2-2ic\tau)}\,\frac{\cosh (2Ms\sqrt{\tau^2-p^2})}{\sqrt{p^2-\tau^2}}\,d\tau
\]
and by the path ${\cal L}$ it is understood that this refers to the steepest descent paths ${\cal C}_1^\pm$ situated in $\Re (\tau)\geq 0$ for the positive and negative components of the integral, respectively. 
When $\alpha=0$, the expression on the right-hand side of (\ref{e32}) agrees with (\ref{e23}) and (\ref{e23a}).

In the first integral in (\ref{e32}) we expand the trigonometric functions to obtain
\begin{eqnarray*}
I_1&=&\sum_{k=0}^\infty \sum_{m=0}^\infty (-)^{m+k}\,\frac{(2Mc)^{2k+1}}{(2k+1)!}\,\frac{(2Ms)^{2m}}{(2m)!} \int_0^p e^{-M\tau^2}\tau^{2k+1} (p^2-\tau^2)^{m-\fr}d\tau\\
&=&\frac{\pi e^{-\rho}}{2}\sum_{k=0}^\infty \sum_{m=0}^\infty \frac{(-)^{m+k}(Mcp)^{2k+1}(Msp)^{2m}}{\g(k+\f{3}{2})\,m!\,\g(m+k+\f{3}{2})}\,{}_1F_1(m+\fs; m+k+\f{3}{2};\rho)
\end{eqnarray*}
upon use of (\ref{e24}) with $\mu=m+\fs$ and application of the duplication formula for the gamma function and Kummer's transformation. Then substitution of the series expansion of the hypergeometric function yields the convergent expansion
\bee\label{e33}
I_1=\frac{\pi e^{-\rho}}{2}\sum_{r=0}^\infty\frac{\rho^r}{r!} \sum_{m=0}^\infty \frac{(-)^m (m+\fs)_r}{m!} (\fs xs)^{2m}\,{\cal H}_{m+r}(xc),
\ee
where ${\cal H}_\nu(x)$ is the scaled Struve function defined in (\ref{e2s}). When $\alpha=0$, only the term corresponding to $m=0$ makes a contribution and we recover the evaluation given in (\ref{e25}). 

With the change of variable $\xi=u/p$, the second integral in (\ref{e32}) can be written as
\[I_2=\int_0^{\xi_0} e^{\rho\xi^2-xc\xi}\ \frac{\cos (sx\sqrt{1+\xi^2})}{\sqrt{1+\xi^2}}\,d\xi,\qquad \xi_0=\frac{2Mu_0}{x},\]
where the value of $u_0$ as $p\ra 0$ is given in (\ref{e35}).
We follow the same procedure employed in Section 2.2 to find (see Appendix A for details)
\begin{eqnarray}
I_2&=&\sum_{k=0}^{n-1}\frac{\rho^k}{k!}\int_0^{\xi_0} \xi^{2k}e^{-xc\xi}\, \frac{\cos (sx\sqrt{1+\xi^2})}{\sqrt{1+\xi^2}}\,d\xi+O(M^{-n})\nonumber\\
&=&\frac{\pi}{2}\bl\{\sum_{k=0}^{n-1}\frac{M^{-k}}{2^{2k} k!}\,C_k(x,\alpha)+O(M^{-n})\br\}+O( e^{-Mc^2})\label{e34}
\end{eqnarray}
when $n< Mc^2$, where
\begin{eqnarray}
C_k(x,\alpha)&=&\frac{2}{\pi}x^{2k}\!\!\int_0^\infty \xi^{2k} e^{-xc\xi}\, \frac{\cos (sx\sqrt{1+\xi^2})}{\sqrt{1+\xi^2}}\,d\xi \\
&=&\frac{2}{\pi} \int_0^\infty t^{2k} e^{-ct} \,\frac{\cos (s\sqrt{x^2+t^2})}{\sqrt{x^2+t^2}}\,dt ;\label{e34a}
\end{eqnarray}
compare (\ref{e29}) in the case $\alpha=0$.

Finally, we determine the contribution $I_s$. The positive component of the integral $I_s$ passes along the steepest descent path through the saddle $\tau_s^+$, where for $\alpha$ bounded away from zero we obtain by
application of the saddle-point method the estimate
\[\bl(\frac{2\pi}{-M\psi''(\tau_s^+)}\br)^{\!\fr} e^{M\psi(\tau_s^+)+\fr i\alpha},\]
where $\psi''(\tau_s^+)\simeq -2$.
From (\ref{e35e}), we have
\[\psi(\tau_s^+)\simeq -e^{-i\alpha}+ip^2 se^{i\alpha/2} \qquad (p\ra 0).\]
The contribution from the negative component is more recessive (being of $O(M^{-1}e^{-Mc^2})$), since it involves a portion of the steepest descent path ${\cal C}_1^-$ {\it not\/} including the saddle point $\tau_s^-$. Hence we find the estimate for large $M$
\bee\label{e38}
e^{\fr\rho} I_s \simeq (\pi/M)^\fr \,e^{-(M-\fr\rho) \cos \alpha}\,\sin\,[(M+\fs\rho) \sin \alpha+\fs\alpha]
\ee
when $0<\alpha\leq\fs\pi$, which agrees with that obtained by Ursell stated in (\ref{e13}). This estimate is exponentially small except when $\alpha$ is near $\fs\pi$. It is also important to stress that the above estimate is not valid as $\alpha\ra 0$; see (\ref{e212}) when $\alpha=0$.

From (\ref{e32}), (\ref{e33}), (\ref{e34}) and (\ref{e38}) we obtain our desired result when $0<\alpha\leq\fs\pi$ given in the following theorem.
\begin{theorem}$\!\!\!.$\ \ Let $x$ and $\rho$ be positive and bounded and let $0<\alpha\leq\fs\pi$, with $s\equiv \sin \fs\alpha$, $c\equiv \cos \fs\alpha$. Then, when $M\equiv x^2/(4\rho)$ is large, we have
\begin{eqnarray}
F(x,\rho,\alpha)&=&-\pi e^{-\fr\rho} \sum_{r=0}^\infty\frac{\rho^r}{r!} \sum_{m=0}^\infty \frac{(-)^m (m+\fs)_r}{m!} (\fs xs)^{2m}\,{\cal H}_{m+r}(xc) \nonumber\\
&&+\pi e^{\fr\rho} \bl\{\sum_{k=0}^{n-1}\frac{M^{-k}}{2^{2k} k!}\,C_k(x,\alpha)+O(M^{-n})\br\}\nonumber\\
&&+ O(e^{-Mc^2})+O(M^{-\fr} e^{-M\cos \alpha})\label{e310}
\end{eqnarray}
where $n$ ($<Mc^2$) is a fixed positive integer, ${\cal H}_\nu(x)$ is the scaled Struve function defined in (\ref{e2s}) and the coefficients $C_k(x,\alpha)$ are defined in (\ref{e34a}).
\end{theorem}

Some examples that demonstrate the precision of the asymptotic formulas in Theorems 1 and 2 will now be given. In Fig.~3 we show the behaviour of the first few coefficients $C_k(x,\alpha)$ as a function of $x$ for a given $\alpha$. Denoting the convergent sum involving Struve functions in (\ref{e310}) by $S_1$ for brevity, then we define
\[
{\cal F}:=F(x,\rho,\alpha)+\pi e^{-\fr\rho} S_1-\pi e^{\fr\rho} \sum_{k=0}^{n-1}\frac{M^{-k}}{2^{2k} k!}\,C_k(x,\alpha).
\]
\begin{figure}[t]
	\begin{center}\includegraphics[width=0.4\textwidth]{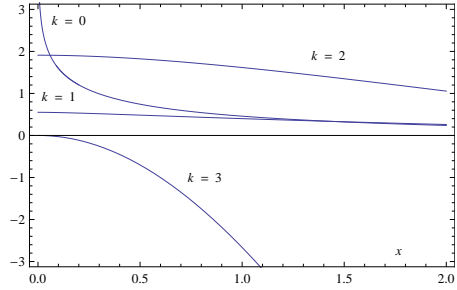}\\
\caption{\small{Behaviour of the first few coefficients $C_k(x,\alpha)$ as a function of $x$ when $\alpha=\pi/6$.}}
	\end{center}
\end{figure}
Table 1 shows the absolute error in the computation of ${\cal F}$ from the expansions in (\ref{e213}) and (\ref{e310}) for different values of $\alpha$ together with the truncation index $n$ employed (satisfying $n<Mc^2$). The convergent sum $S_1$ was computed with sufficient terms commensurate with the overall accuracy of the computation. It is found that these errors are comparable to those resulting from optimal truncation of the Bessho expansion in (\ref{e12}) (without the integral term). 

\begin{table}[th]
\caption{\footnotesize{Values of the absolute error in the computation of ${\cal F}$ and the truncation index $n$.}}
\begin{center}
\begin{tabular}{|l|lc|lr|}
\hline
&&&&\\[-0.3cm]
\mcol{1}{|c|}{} & \mcol{2}{|c|}{$x=0.40,\ \rho=0.005$} & \mcol{2}{c|}{$x=1.00,\ \rho=0.020$}\\  
\mcol{1}{|c|}{} & \mcol{2}{|c|}{($M=8$)} & \mcol{2}{|c|}{($M=12.5$)}\\
\mcol{1}{|c|}{$\alpha/\pi$}& \mcol{1}{c}{$|{\cal F}|$} & \mcol{1}{c|}{$n$} & \mcol{1}{c}{$|{\cal F}|$} & \mcol{1}{c|}{$n$} \\
[.1cm]\hline
&&&&\\[-0.3cm]
0 & $6.368\times 10^{-6}$ & 8 & $2.613\times 10^{-7}$ & 12\\
0.10 & $3.146\times 10^{-6}$ & 5 & $1.998\times 10^{-6}$ & 12\\
0.20 & $3.146\times 10^{-6}$ & 6 & $1.899\times 10^{-5}$ & 11\\
0.25 & $4.687\times 10^{-3}$ & 5 & $1.428\times 10^{-5}$ & 9\\
0.30 & $2.976\times 10^{-3}$ & 3 & $2.890\times 10^{-4}$ & 9\\
0.40 & $4.326\times 10^{-2}$ & 1 & $7.928\times 10^{-4}$ & 8\\
[.2cm]\hline
\end{tabular}
\end{center}
\end{table}
\vspace{0.6cm}

\begin{center}
{\bf Appendix A: \ Derivation of bounds on the contribution from the imaginary axis for $0\leq\alpha\leq\fs\pi$}
\end{center}
\setcounter{section}{1}
\setcounter{equation}{0}
\renewcommand{\theequation}{\Alph{section}.\arabic{equation}}
The contribution to $F(x,\rho,\alpha)$ arising from the imaginary $\tau$-axis over the interval $[0, iu_0]$ is
expressed in the form
\bee\label{b1}
\sum_{k=0}^{n-1}\frac{\rho^k}{k!} \bl(\int_0^\infty-\int_{\xi_0}^\infty\br) \xi^{2k}e^{-xc\xi}\,\frac{\cos (sx\sqrt{1+\xi^2})}{\sqrt{1+\xi^2}}\,d\xi+R_n,
\ee
where
\[R_n=\frac{\rho^n}{n!}\int_0^{\xi_0} \xi^{2n} e^{-xc\xi+\theta\rho\xi^2}\,\frac{\cos (sx\sqrt{1+\xi^2})}{\sqrt{1+\xi^2}}\,d\xi, \quad 0<\theta<1\]
and $\xi_0=2Mu_0/x$. This results from expansion
of $e^{\rho\xi^2}$ as a finite series of $n$ terms together with a remainder in Lagrange's form; see (\ref{e2L}). The quantity $u_0$
is the intercept on the imaginary axis of the steepest descent path through the saddle labelled $\tau_s^+$.
From (\ref{e35}) we have $u_0=c(1-\fs p^2\tan^2 \fs\alpha)$ as $p\ra 0$, so that $u_0\leq c$. When $\alpha=0$, we have $u_0=1$, $\xi_0=2Mc/x$ and the above expressions reduce to those given in (\ref{e26}) and (\ref{e26a}).

The remainder term $R_n$ satisfies
\[|R_n|<\frac{\rho^n}{n!}\int_0^{\xi_0}\xi^{2n-1}e^{-xc\xi+\theta\rho\xi^2}d\xi=\frac{\rho^n}{n!}\int_0^{\xi_0}\xi^{2n-1}e^{-xc\xi/2}\cdot g(\xi)\,d\xi,\]
where $g(\xi)=\exp [-\fs xc\xi+\theta\rho\xi^2]$. It is easily shown that $g(\xi)$ is monotonically decreasing on $[0, \xi_0]$ when $0<\theta\leq\theta_0$, $\theta_0\equiv c/(2u_0)$, and that there is a single minimum at $\xi=Mc/(\theta x)$ when $\theta_0\leq\theta<1$. Hence it follows that $g(\xi)\leq 1$ on the interval $[0,\xi_0]$ since, with $u_0\leq c$,
\[g(\xi_0)=e^{-Mu_0^2((c/u_0)-\theta)}<e^{-Mu_0^2(1-\theta)}<1.\]
%\[g(\xi_0)=e^{-Mu_0^2((c/u_0)-\theta)}<e^{-Mu_0^2(1-\theta)}\]
%\[g^*\leq g(\xi)\leq 1\quad (0\leq\xi\leq\xi_0),\]
%where, with $\theta_0\equiv c/(2u_0)$,
%\[g^*=\bl\{\begin{array}{ll}\!\!g(\xi_0)=e^{-Mu_0^2((c/u_0)-\theta)}<e^{-Mu_0^2(1-\theta)}& (0<\theta\leq \theta_0)\\ \!\!g(\xi_*)=e^{-Mc^2/(4\theta)} & (\theta_0\leq\theta<1)\end{array}\]
Then we have the bound
\bee\label{e27}
|R_n|<\frac{\rho^n}{n!} \int_0^{\xi_0} \xi^{2n-1} e^{-xc\xi/2}d\xi< \frac{M^{-n} \g(2n)}{n!}
\ee
upon replacing the upper limit of integration by $\infty$.

The tail $T$ of the finite sum in (\ref{b1}) involving the integral over $[\xi_0,\infty)$ can be estimated as follows:
\begin{eqnarray*}
|T|&<&\sum_{k=0}^{n-1}\frac{\rho^k}{k!} \int_{\xi_0}^\infty \frac{\xi^{2k} e^{-xc\xi}}{\sqrt{1+\xi^2}}\,d\xi<\sum_{k=0}^{n-1}\frac{\rho^k}{k!}\int_{\xi_0}^\infty \xi^{2k-1} e^{-xc\xi}d\xi\\
&=&\sum_{k=0}^{n-1} \frac{(Mc^2)^{-k}}{2^{2k} k!} \int_{2Mu_0c}^\infty w^{2k-1} e^{-w}\,dw=\sum_{k=0}^{n-1} \frac{(Mc^2)^{-k}}{2^{2k} k!}\,\g(2k,2Mu_0c),
\end{eqnarray*}
where $\g(a,\chi)$ denotes the upper incomplete gamma function. We now suppose that $n<Mu_0c\leq Mc^2$ and employ the simple bound (see Appendix B) for $a>0$, $\chi>0$ given by $\g(a,\chi)\leq 2\chi^a e^{-\chi}$ when $0\leq a\leq \chi$. This yields the bound, for $n<Mu_0c\leq Mc^2$,
\begin{eqnarray}
|T|&<&2e^{-2Mu_0c} \sum_{k=0}^{n-1} \frac{(Mu_0^2)^k}{k!}<2e^{-2Mu_0c+Mu_0^2}\nonumber\\
&<&2e^{-Mu_0c}=O(e^{-Mc^2})\qquad(n<Mc^2).\label{b2}
\end{eqnarray}
If $n$ is finite, the above bound can be replaced by 
\bee\label{a4}
|T|< 2e^{-2Mu_0c} \sum_{k=0}^{n-1} (Mu_0^2)^k=O(M^ne^{-2Mc^2}).
\ee

\vspace{0.6cm}
\begin{center}
{\bf Appendix B: \ Derivation of a bound on $\g(a,\chi)$}
\end{center}
\setcounter{section}{2}
\setcounter{equation}{0}
\renewcommand{\theequation}{\Alph{section}.\arabic{equation}}
The derivation of this bound follows closely that given in \cite{U2}; see also \cite[p.~76]{P}. Suppose that $a>0$ and $\chi>0$. Then \cite[p.~174]{DLMF}
\[\g(a+1,\chi)=\int_\chi^\infty e^{-t}t^a dt=e^{-\chi} \chi^{a+1} \int_0^\infty e^{-u\chi} (1+u)^a du.\]
When $-1\leq a\leq 0$, we therefore have
\[\g(a+1,\chi)\leq e^{-\chi} \chi^{a+1}\int_0^\infty e^{-u\chi}du\leq e^{\chi} \chi^{a+1} \quad (\chi\geq 1).\]
When $0\leq a\leq 1$, we have $(1+u)^a\leq 1+u$ for $u\geq 0$, so that
\begin{eqnarray*}
\g(a+1,\chi)&\leq & e^{-\chi} \chi^{a+1}\int_0^\infty e^{-u\chi}(1+u)\,du\leq e^{\chi} \chi^{a+1}\bl(\frac{1}{\chi}+\frac{1}{\chi^2} \br)\\
&\leq & 2e^{-\chi} \chi^{a+1}\qquad (\chi\geq 1).
\end{eqnarray*}
When $a\geq 1$, we have
\begin{eqnarray*}
\g(a+1,\chi)&\leq & e^{-\chi} \chi^{a+1}\int_0^\infty \{e^{-u}(1+u)\}^a\,du\\
&\leq& e^{-\chi} \chi^{a+1} \int_0^\infty e^{-u}(1+u) du=2e^{-\chi} \chi^{a+1},
\end{eqnarray*}
the first inequality holding when $a\leq\chi$ and the second inequality when $a\geq 1$, since the integrand $e^{-u}(1+u)\leq 1$ on $[0,\infty)$. The resulting bound therefore holds for $1\leq a\leq\chi$.

Collecting together these results and replacing $a$ by $a-1$, we therefore obtain the upper bound
\bee\label{c1}
\g(a,\chi)\leq 2e^{-\chi} \chi^{a} \qquad (0\leq a\leq \chi)
\ee
provided $\chi\geq 1$.

\end{document}